\documentclass[12pt, english]{article}

\usepackage{amsmath}
\usepackage{graphicx}
\usepackage{amssymb}
\usepackage{subcaption}
\usepackage{dsfont}
\usepackage{url}
\newcommand{\beqn}{\begin{eqnarray*}}
\newcommand{\eeqn}{\end{eqnarray*}}
\newcommand{\bneqn}{\begin{eqnarray}}
\newcommand{\eneqn}{\end{eqnarray}}
\newcommand{\parens}[1]{\left(#1\right)}

\newcommand{\bracks}[1]{\left[#1\right]}

\newcommand{\abss}[1]{\left|#1\right|}
\newcommand{\expe}[1]{\mathbb{E}\bracks{#1}}
\newcommand{\expesub}[2]{\mathbb{E}_{#1}\bracks{#2}}

\renewcommand{\exp}[1]{\mathrm{exp}\parens{#1}}
\newcommand{\indic}[1]{\mathds{1}_{#1}}
\renewcommand{\sin}[1]{\text{sin}\parens{#1}}

\title{On occupation times of the first and third \\ quadrants for planar Brownian motion}
\author{Philip A. Ernst \footnote{Department of Statistics, Rice University} and Larry Shepp \footnote{Deceased April 23, 2012}}
\date{\today}

\begin{document}
\maketitle

\begin{abstract}

In 1988, Bingham and Doney in 1988 (\cite{Bingham}) presented the applied probability community with a question which is very simply stated, yet is extremely difficult to solve: what is the distribution of the quadrant occupation time of planar Brownian motion? In this paper, we study an alternate formulation of this long-standing open problem: let $X(t), Y(t), t \geq 0$ be standard Brownian motions starting at $x,y$ respectively. Find the distribution of the total time $T=Leb\{t \in [0,1]: X(t) \times Y(t) >0\}$, when $x=y=0$, i.e., the occupation time of the union of the first and third quadrants. If two adjacent quadrants are used, the problem becomes much easier and the distribution of $T$ follows the arcsine law.  
\end{abstract}

\section{Introduction}
In this paper, we consider a long-standing open problem in the applied probability literature: what is the quadrant occupation time of planar Brownian motion? This question had intrigued Larry Shepp since 1995 (see \cite{SheppSem}). Formally, let $T$ be the total time that the vector process $X(t) = (W_1(t),W_2(t))$ on $0 \leq t \leq 1$ is in the first quadrant; the task is to find the distribution of $T$. In 1988, Bingham and Doney remarked on p. 121 of \cite{Bingham} that ``in no case to our knowledge is the law of $T$ known explicitly.'' Using independence of coordinate processes, the authors obtained the first two moments of $T$ and provided a solution for the third moment, the latter of which was corrected by \cite{Desbois}. The author of \cite{Desbois} generalized the aforementioned quadrant problem by considering the occupation time spent in a wedge of apex $O$ and angle $\theta$. Analytical results for general $\theta$ were provided for both the second and third moments, and the fourth moment for the quadrant problem ($\theta=\pi/2$) was obtained. \\
\indent Despite these new additions to the literature, the author of \cite{Desbois} concludes that ``our feeling is that the occupation time problem for Brownian motion is far from being understood as soon as we leave one-dimensional or quasi-one-dimensional (graphs) situations. Clearly, new ideas are needed if we want to tackle this problem.'' Our work, through its use of Kontorovich-Lebedev transforms and pasting of solutions, offers a new, but ultimately incomplete, approach on this long-standing problem. For references on the Kontorovich-Lebedev transform, we refer to the reader to the following sources: \cite{KL2}, \cite{KL1}, \cite{KL4}, and \cite{KL3}.
\section{Main Results}

\subsection{Setup}
We consider the following alternative formulation of Bingham and Doney's (\cite{Bingham}) quadrant occupation of planar Brownian motion problem. Let
\beqn
X(t), Y(t), t \geq 0
\eeqn
be standard Brownian motions starting at $x,y$ respectively. We wish to find the distribution of the total time $T=Leb\{t \in [0,1]: X(t) \times Y(t) >0\}$, when $x=y=0$, i.e., the occupation time of the union of the first and third quadrants. If two adjacent quadrants are used, the problem becomes much easier and the distribution follows the arcsine law (\cite{Levy}).

The Feynman-Kac theorem states (see \cite{Karatzas}) that

\bneqn
U(x,y)=\expesub{x,y}{\int_0^\infty \exp{-\alpha t-\lambda\int_0^t \indic{\parens{X(u) \times Y(u)>0}}du}dt}
\eneqn
is the bounded solution of the Helmholtz partial differential equation in each quadrant,

\bneqn
\frac{1}{2}(U_{xx}+U_{yy})(x,y)-\beta(x,y)U(x,y)+1 \equiv 0,
\eneqn
where $\beta(x,y)=\beta_1=\alpha+ \lambda$ for $(x,y) \in Q_1$, $Q_3$ (the first and third quadrants, respectively) and $\beta(x,y)=\beta_2=\alpha, (x,y) \in Q_2, Q_4$ (the second and fourth quadrants, respectively). The function $U$ must be twice differentiable interior to each quadrant, continuously differentiable overall, and uniformly bounded. If we can find $U$, then we know $U(0,0)$, and then we will have:

\bneqn
U(0,0)= \expe{\int_0^\infty e^{-\alpha t-\lambda tT}dt}=\expe{\frac{1}{\alpha+\lambda T}}.
\eneqn

We turn to finding the Kontorovich-Lebedev solution in each quadrant. Let $x=r\cos\theta$, $y=r\sin\theta$, and set $V(r,\theta)=U(x,y)$. As is well known,

\bneqn
U_{xx}+U_{yy}=V_{rr}+\frac{1}{r}V_r+ \frac{1}{r^2}V_{\theta\theta}.
\eneqn
The modified Bessel function $v(r)=\kappa_{iv}(r)$ satisfies the ordinary differential equation (see \cite{Ober}):
\bneqn
r^2v^{''}(r)+rv'(r)-(r^2-\nu^2)v(r)=0.
\eneqn
It now can be easily checked that, for any functions $f(\nu), g(\nu)$,

\bneqn
V(r,\theta)=\frac{1}{\beta}+\int_0^\infty f(\nu)\kappa_{iv}(r\sqrt{2\beta})\sinh(\nu\theta)d\nu+\int_0^\infty g(\nu)\kappa_{iv}(r\sqrt{2\beta})\cosh(\nu\theta)d\nu
\eneqn
solves the differential equation:
\bneqn
\frac{1}{2}\parens{V_{rr}(r,\theta)+\frac{1}{r}V_r(r,\theta)+\frac{1}{r^2}V(r,\theta)}-\beta V(r, \theta)+1=0,
\eneqn
with a different choice of $f,g, \beta$ in each quadrant, which we must paste together to satisfy the needed smoothness. For convenience, we can use any two linearly independent combinations of $\sinh, \cosh$, etc. We now proceed to do so.

\subsection{Pasting}
Our strategy is to use ``pasting'' of the solutions. We denote $f_i$ and $g_i$ as the densities for each of the two linear combinations of $\sinh, \cosh$, respectively, in the $i$th quadrant. In:

\beqn
Q_1&=& \{(r,\theta): r >0, 0<\theta<\pi/2\}
\eeqn
we set:
\small
\bneqn
V(r,\theta)=\frac{1}{\beta_1}+\int_0^\infty f_1(\nu)\kappa_{iv}\parens{r\sqrt{2\beta_1}}\sinh\parens{\nu\parens{\frac{\pi}{2}-\theta}}d\nu+\int_0^\infty g_1(\nu)\kappa_{iv}\parens{r\sqrt{2\beta_1}}\sinh (v\theta)d\nu.
\eneqn
\normalsize
In

\bneqn
Q_2&=& \{(r,\theta): r >0, \frac{\pi}{2}<\theta<\pi\}
\eneqn
we set:

\footnotesize

\bneqn
V(r,\theta)=\frac{1}{\beta_2}+\int_0^\infty f_2(\nu)\kappa_{iv}\parens{r\sqrt{2\beta_2}}\sinh\parens{\nu\parens{\pi-\theta}}d\nu+\int_0^\infty g_2(\nu)\kappa_{iv}\parens{r\sqrt{2\beta_2}}\sinh \parens{v\parens{\theta-\frac{\pi}{2}}}d\nu.
\eneqn
\normalsize
By symmetry, $U(x,y)=U(y,x)=U(-x,-y)$, and so $g_j \equiv f_j$ for $j=1,2,3,4$, $f_3=f_1,f_4=f_2$.
\subsection{Consequences of Continuity and Continuous Differentiability on the Axes}
Note that:
\beqn
V\parens{r,\frac{\pi}{2}-0}&=&\frac{1}{\beta_1}+\int_0^\infty f_1(\nu)\sinh\parens{\frac{\nu\pi}{2}}\kappa_{iv}(r\sqrt{2\beta_1})d\nu,\\
V\parens{r,\frac{\pi}{2}+0}&=&\frac{1}{\beta_2}+\int_0^\infty f_2(\nu)\sinh\parens{\frac{\nu\pi}{2}}\kappa_{iv}(r\sqrt{2\beta_2})d\nu,
\eeqn
and the right-hand sides of these equations are equal. The derivatives, taken with respect to $\theta$, are:
\beqn
V_{\theta}\parens{r,\frac{\pi}{2}-0}&=&\int_0^\infty f_1(\nu)\nu\parens{\cosh\parens{\frac{\nu\pi}{2}}-1}\kappa_{iv}(r\sqrt{\beta_1})d\nu,\\
V_{\theta}\parens{r,\frac{\pi}{2}+0}&=&\int_0^\infty f_2(\nu)\nu\parens{-\cosh\parens{\frac{\nu\pi}{2}}+1}\kappa_{iv}(r\sqrt{\beta_2})d\nu,
\eeqn
and the right-hand sides of these equations are equal.

\subsubsection{Solving the Above Smoothness Equations}
We assume that there are signed measures, $\mu_j(dz),\, j=1,2,$ such that:

\bneqn \label{cooleqn}
f_j(\nu)=\frac{2}{\pi}\frac{\cosh(\nu\frac{\pi}{2})}{\sinh(\nu\frac{\pi}{2})}\parens{\frac{-1}{\beta_j}}+\int_0^\infty \mu_j\, \sin{\nu z}dz.
\eneqn
We then can plug in $f_j$ into the first two equations expressing continuity on the $y$-axis. The first term in $f_j$ kills the term $\frac{1}{\beta_j}$ because the Kontorovich-Lebedev transform of $\cosh(\nu\frac{\pi}{2})$ is, for every $y$,

\bneqn
\frac{2}{\pi}\int_0^\infty\cosh\parens{\nu\frac{\pi}{2}}\kappa_{iv}(y)d\nu=e^{-y\cos{\frac{\pi}{2}}}\equiv 1,
\eneqn
(see \cite{Ober} p.242). This leaves the following equation for the sine transforms of $\mu_j$ for $V(r,\theta)$ to be continuous at each $r$ when $\theta=\frac{\pi}{2}$:

\bneqn
V\parens{r,\frac{\pi}{2}}&=&\int_0^\infty \mu_1dz\int_0^\infty\sin{\nu z}\sinh\parens{\frac{\nu \pi}{2}}\kappa_{iv}(r\sqrt{2\beta_1})d\nu\\
&=&\int_0^\infty \mu_2dz\int_0^\infty\sin{\nu z}\sinh\parens{\frac{\nu \pi}{2}}\kappa_{iv}(r\sqrt{2\beta_2})d\nu.
\eneqn
Identity (8) on p.244 of \cite{Ober} states that

\beqn
\int_0^\infty\sin{\nu z}\sinh\parens{\frac{\nu\pi}{2}}\kappa_{iv}(y)d\nu=\frac{\pi}{2}\sin{y\sinh(z)}.
\eeqn
This gives us the continuity equation linking $\mu_1, \mu_2$ as follows:

\bneqn
\int_0^\infty \sin{r\sqrt{2\beta_1}\sinh(z)}\mu_1dz=\int_0^\infty \sin{r\sqrt{2\beta_2}\sinh(z)}\mu_2dz.
\eneqn
We now define the change of variables $z'(z)=\phi(z), z\geq 0$ such that:

\bneqn
\sqrt{2\beta_1}\sinh(z)=\sqrt{2\beta_2}\sinh(\phi(z)).
\eneqn
Explicitly,
\bneqn
\phi(z)=\log\parens{{\sqrt{\frac{\beta_1}{\beta_2}}\sinh z+\sqrt{\frac{\beta_1}{\beta_2}\sinh^2(z)+1}}}.
\eneqn
We then have, for every $r\geq 0$,

\bneqn
\int_0^\infty \sin{r\sqrt{2\beta_1}\sinh(z)}\mu_1dz=\int_0^\infty \sin{r\sqrt{2\beta_1}\sinh(z)}\mu_2d\phi(z).
\eneqn
Since the sine transform is unique on a half interval, we have:

\bneqn
\mu_1(dz)=\mu_2(d\phi(z)), z\geq 0.
\eneqn
Thus, we have expressed one relationship between $\mu_1$ and $\mu_2$. We now need a second equation linking $\mu_1$ and $\mu_2$, and we use the equation obtained by using the continuity of the derivative on the positive $y$-axis. 

\subsubsection{Continuity of the Derivative of $V$ on $\theta$ at $\theta=\frac{\pi}{2}$}
Noting that:

\beqn
V_{\theta}\parens{r,\frac{\pi}{2}-0}&=&\int_0^\infty f_1(\nu)\nu\parens{\cosh\parens{\frac{\nu\pi}{2}}-1}\kappa_{iv}(r\sqrt{2\beta_1})d\nu,\\
V_{\theta}\parens{r,\frac{\pi}{2}+0}&=&-\int_0^\infty f_2(\nu)\nu\parens{\cosh\parens{\frac{\nu\pi}{2}}-1}\kappa_{iv}(r\sqrt{2\beta_2})d\nu,
\eeqn
we immediately see that the right-hand sides are equal for all $r \geq 0$ (the derivative of $V$ on $\theta$ is continuous at $\theta=\frac{\pi}{2}$). Placing the expressions for the $f_j$, given in (\ref{cooleqn}), in terms of the unknown $\mu_j$, into the right-hand side of each of the above equations gives:

\bneqn
0=\sum_{j=1}^2\int_0^\infty\bracks{\frac{2}{\pi}\frac{\cosh(\nu\frac{\pi}{2})}{\sinh(\nu\frac{\pi}{2})}\parens{\frac{-1}{\beta_j}}+\int_0^\infty \mu_j\sin{\nu z}dz}\nu\parens{\cosh\parens{\frac{\nu\pi}{2}}-1}\kappa_{iv}(r\sqrt{2\beta_j})d\nu.
\eneqn
From \cite{Ober}, p.244, equation (7), for $a$ real or $\abss{a}\leq \frac{\pi}{2}$, 

\bneqn
\frac{2}{\pi}\int_0^\infty\nu \sin{a\nu}\kappa_{i\nu}(y)d\nu=ye^{-y\cosh a}\sinh a.
\eneqn
We then use the identities:\\
\bneqn 
\frac{\cosh\parens{\frac{\nu \pi}{2}}}{\sinh\parens{\frac{\nu \pi}{2}}}\parens{\cosh\parens{\frac{\nu \pi}{2}}-1}&=&\sinh\parens{\frac{\nu\pi}{2}}-\tanh\parens{\frac{\nu\pi}{4}}, \nonumber \\ 
\sin{\nu z}\parens{\cosh\parens{\frac{\nu \pi}{2}}-1}&=&\frac{1}{2}\sin{\nu\parens{z+\frac{i\pi}{2}}}+\frac{1}{2}\sin{\nu\parens{z-\frac{i\pi}{2}}}-\sin{\nu z}, \nonumber
\eneqn
and link $\mu_1$ and $\mu_2$ as follows:

\scriptsize
\bneqn \label{cooldude}
0&=&\sum_{j=1}^2\int_0^\infty\frac{2}{\pi}\parens{\frac{-1}{\beta_j}}\parens{\sinh\parens{\frac{\nu\pi}{2}}-\tanh\parens{\frac{\nu\pi}{4}}}\nu\kappa_{iv}(r\sqrt{2\beta_j})d\nu \nonumber \\
&+&\sum_{j=1}^2\int_0^\infty \mu_jdz \int_0^\infty \bracks{\frac{1}{2}\sin{\nu\parens{z+\frac{i\pi}{2}}}+\frac{1}{2}\sin{\nu\parens{z-\frac{i\pi}{2}}}-\sin{\nu z}\nu\kappa_{iv}(r\sqrt{2\beta_j})}d\nu. 
\eneqn
\normalsize
Since $\sinh(z+2i\sqrt{\pi})=i\cosh z$ and $\cosh(z+i\pi/2)=-i\sinh z$, we obtain from (\ref{cooldude}), with $a=i\frac{\pi}{2}$, $a=z\pm i\frac{\pi}{2}$, respectively:

\normalsize

\bneqn
0&=&\sum_{j=1}^2\parens{-r\sqrt{\frac{2}{\beta_j}}}+\frac{2}{\pi \beta_j}\int_0^\infty \nu \tanh \parens{\frac{\nu\pi}{4}}\kappa_{iv}(r\sqrt{2\beta_j})d\nu \nonumber\\
&+&\sum_{j=1}^2\frac{\pi}{2}\int_0^\infty \mu_j\frac{r\sqrt{2\beta_j}}{2}\bracks{e^{r\sqrt{2\beta_j}i \sinh z}(i \cosh z)+e^{-r\sqrt{2\beta_j}i \sinh z}(-i \cosh z)-2e^{-r\sqrt{2\beta_j}\cosh z}\sinh z}dz. \nonumber
\eneqn
\normalsize
We now use the substitution $z'(z)=\phi(z)$, implicitly defined by 

\beqn
\sqrt{2\beta_1}\sinh z= \sqrt{2\beta_2}\sinh z'
\eeqn
in the second integral of the last display, with $j=2$. Combined with the fact that 

\beqn
\mu_2(d\phi(z))=\mu_1(dz)
\eeqn
we obtain:
\small

\bneqn \label{lastone}
&&\frac{\pi}{2}\int_0^\infty \mu_1 r\sqrt{2\beta_1}\bracks{\sin{r\sqrt{2\beta_1}\sinh z}\cosh z + e^{-r\sqrt{2\beta_1}\cosh z}\sinh z}dz \nonumber \\
&+& \frac{\pi}{2}\int_0^\infty \mu_1 r\sqrt{2\beta_2} \bracks{\sin{r\sqrt{2\beta_1}\sinh z} \cosh \phi (z)+ e^{-r\sqrt{2\beta_2}\cosh \phi (z)}\sinh \phi(z)}dz \nonumber \\ 
&=& \sum_{j=1}^2\parens{-r\sqrt{\frac{2}{\beta_j}}}+\frac{2}{\pi \beta_j}\int_0^\infty\tanh \parens{\frac{\nu\pi}{4}}\nu\kappa_{iv}(r\sqrt{2\beta_j})d\nu.
\eneqn

\normalsize
Despite much joint effort, further explicit calculations beyond Equation (\ref{lastone}) quickly become intractable. It is our opinion that an explicit solution to this long-standing problem may not be possible. Nonetheless, we have successfully reduced the problem to that where an analyst of special functions could pick up where we have left off. The problem now becomes one of finding the relationship between functions $f$ and $g$ if their Kontorovich-Lebedev transforms $F$ and $G$ satisfy $F(r) = G(cr)$ for all $r$ with $c$ given.

\section{Remarks}
\begin{enumerate}
\item Professor Terry Lyons is credited (Personal communication with Professor Nick Bingham) with saying that the simplest case beyond the half-plane (which reduces to the arc-sine law in one dimension) is the third plane $0 < \theta < 2 \pi/3$.  This could help one to compare the differential equations we obtain in (\ref{lastone}). 
\item One can consider the random occupation measure generated on the unit circle (or a sphere in higher dimensions) by the angular part of a Brownian motion starting at 0 and running for time 1. Some results regarding this random measure, and its relation to the angle of the Brownian motion at time 1, were obtained in \cite{Pemantle}.
\item Another natural problem in two dimensions is to find the law of the occupation time $A_u$ of an interval of length $2 \pi u$ around the unit circle, for $0 < u < 1.$  Some problems involving cyclically stationary local-time processes were treated in \cite{Pitman}.
\end{enumerate}

\section{Acknowledgments}
We thank Professor Nick Bingham, Professor Jim Pitman, and Professor Robin Pemantle for helpful discussions. We are extremely appreciative of the work of an anonymous referee, whose very helpful report greatly strengthened the quality of this work.

\bibliographystyle{plain}
\bibliography{VF}

\end{document}